\documentstyle[amsfonts,12pt]{article}
\input amssym.def
\newcommand{\erarrow}{\stackrel{\sim}{\rightarrow}}

\newcommand{\tka}{{\widetilde{\kappa}}}
\newcommand{\totimes}{{\widetilde{\otimes}}}
\newcommand{\bom}{\overline{\omega}}
\newcommand{\ZNg}{Z_{N_g}}
\newcommand{\ltimes}{{\,\triangleright\!\!\!<}}
\newcommand{\AGP}{\bbC[G]\ltimes A}

\newcommand{\bul}{{\bullet}}
\newcommand{\al}{{\alpha}}
\newcommand{\la}{{\lambda}}
\newcommand{\h}{{\hbar}}

\newcommand{\mm}{{\mathfrak{m}}}

\newcommand{\mb}{{\mathfrak{b}}}

\newcommand{\om}{{\omega}}

\newcommand{\Om}{{\Omega}}

\newcommand{\si}{{\sigma}}
\newcommand{\ga}{{\gamma}}

\newcommand{\ve}{{\varepsilon}}
\newcommand{\ka}{{\kappa}}
\newcommand{\G}{{\Gamma}}
\newcommand{\pa}{{\partial}}

\newcommand{\W}{{\cal W}}

\newcommand{\cD}{{\cal D}}

\newcommand{\cA}{{\cal A}}

\newcommand{\cW}{{\cal W}}
\newcommand{\cO}{{\cal O}}
\newcommand{\cC}{{\cal C}}
\newcommand{\cK}{{\cal K}}

\newcommand{\bcC}{\overline{\cal C}}

\newcommand{\bC}{{\overline{C}}}

\newcommand{\bbC}{{\Bbb C}}

\newcommand{\n}{{\nabla}}

\newcommand{\de}{{\delta}}
\newcommand{\D}{{\Delta}}

\newcommand{\AG}{A[G]}
\newcommand{\bAG}{A[G]^{op}}

\date{}
\newtheorem{defi}{Definition}
\newtheorem{pred}{Proposition}

\newtheorem{teo}{Theorem}
\newtheorem{conjecture}{Conjecture}
\newtheorem{example}{Example}

\begin{document}

\vspace{-10cm}
\begin{flushright}
 \begin{minipage}{1.2in}
ITEP-TH-43/04
 \end{minipage}
\end{flushright}
\vspace{1.1cm}

\begin{center}
{\huge\bf Hochschild cohomology of quantized symplectic orbifolds\\[0.3cm]
and the Chen-Ruan cohomology}\\[0.7cm]
Vasiliy Dolgushev and Pavel Etingof \\[0.5cm]
{\it Department of Mathematics, MIT,} \\
{\it 77 Massachusetts Avenue,} \\
{\it Cambridge, MA, USA 02139-4307,}\\
E-mails: vald@mit.edu, etingof@math.mit.edu\\[1cm]
\end{center}

\begin{abstract}
We prove the additive version of the conjecture
proposed by Ginzburg and Kaledin in \cite{GK}.
This conjecture states that if $X/G$ is an orbifold
modeled on a quotient of a smooth affine symplectic
variety $X$ (over $\bbC$) by a finite group $G\subset Aut(X)$
and $A$ is a $G$-stable quantum algebra of functions on $X$
then the graded vector space $HH^{\bul}(A^G)$
of the Hochschild cohomology of the algebra
$A^G$ of invariants is isomorphic to
the graded vector space $H^{\bul}_{CR}(X/G)((\h))$
of the Chen-Ruan (stringy) cohomology
of the orbifold $X/G$.
\end{abstract}
~\\[0.3cm]
MSC-class: 16E40, 53D55.

\section{Introduction}
A topological space which is obtained as a
quotient of another space with respect to
an action of a group is still a challenging object
in modern geometry. Orbifolds provide us with
an important class of examples of such spaces.

In this paper we consider a smooth affine algebraic variety $X$
(over $\bbC$) endowed with an algebraic symplectic structure
and with an action $\pi$ of a finite group $G$ which preserves
the symplectic structure. Given these data it is not hard to
construct a $G$-invariant (associative) star-product $\star$ in
the vector space $A_+ = \cO(X)[[\h]]$ of formal power
series over regular functions on $X$, which deforms the
commutative multiplication on $\cO(X)$ in the sense of
the deformation quantization \cite{Bayen}, \cite{Ber}.
If $A=A_+[\h^{-1}]$ denotes the localization
of $A_+$ then it is natural to refer to
the algebra $A^G$ of invariants as the algebra of
quantum functions of the orbifold modeled on $X/G$.

The main result of this paper is
\begin{teo}
\label{ONA}
If the group $G$ acts faithfully on $X$
then there exists an isomorphism
\begin{equation}
\label{main}
HH^{\bul}(A^{G})\cong H^{\bul}_{CR}(X/G, \bbC((\h)))
\end{equation}
between the graded vector space $HH^{\bul}(A^G)$ of
Hochschild cohomology of the quantum algebra
$A^G$ of functions on $X/G$ and the graded vector
space $H^{\bul}_{CR}(X/G, \bbC((\h)))$ of Chen-Ruan
(=stringy) cohomology \cite{CR} of the orbifold $X/G$
with coefficients in the field $\bbC((\h))$\,.
\end{teo}
This statement is an additive version of
the conjecture stated in paper \cite{GK}.
More precisely, in \cite{GK} V. Ginzburg and D. Kaledin
conjectured the existence of an isomorphism (\ref{main})
which also respects the corresponding multiplicative
structures.

We have to mention that the algebras we deal with are
defined over the field $\bbC((\h))$ of formal Laurent
power series in an auxiliary variable $\h$.
Therefore our algebras are endowed with the
natural $\h$-adic topology and, in the definition
of Hochschild homological complex, we have to replace
the algebraic tensor product by the one completed in
this topology.

We would like to mention recent
results of N. Ganter who proposed in her PhD
thesis \cite{Nora} a homotopic refinement of the
theory of the orbifold elliptic genus and the
Dijkgraaf-Moore-Verlinde-Verlinde (DMVV)
product formulas \cite{DMVV}. In her thesis
she also defines higher chromatic versions of
the orbifold elliptic genus and proves the DMVV-type
formula for these versions of the
orbifold elliptic genus. In recent paper \cite{chiral}
E. Frenkel and M. Szczesny proved that the graded
vector space of the Chen-Ruan orbifold cohomology
is isomorphic to the graded vector space of the
cohomology of the chiral De Rham complex,
which can be naturally associated with an orbifold
modeled on a quotient of a smooth complex variety.
In \cite{chiral} the
authors also proposed a formula which expresses
the orbifold elliptic genus in terms of a
partition function of the corresponding
chiral vertex algebra. It would be interesting
to generalize this formula  of E. Frenkel
and M. Szczesny to the higher chromatic versions of the
orbifold elliptic genus proposed in \cite{Nora}
by N. Ganter.

The authors of this paper were unaware about the
results of X. Tang who computed in his thesis \cite{Tang}
Hochschild and periodic cyclic homology of the quantum algebra of
functions on a proper \'etale groupoid equipped with an invariant symplectic
structure. Combining the
result of X. Tang (see theorem $4.3.2$) together with the Van Den Bergh duality
\cite{VB} one can easily deduce theorem \ref{ONA}.
However, we have to mention that, in order to
compute Hochschild homology of the above quantum algebra,
X. Tang used the formality theorem for Hochschild
chains \cite{chains} and we do use any version of
the formality theorem in our proof.

The paper is organized as follows.
In the following section we outline a
sketch of the proof of theorem \ref{ONA}.
This proof is based on three propositions
\ref{VdB}, \ref{FLT}, and \ref{ONO}. The
first proposition is the Van Den Bergh duality
\cite{VB} between Hochschild cohomology
and Hochschild homology. The second proposition
is the decomposition theorem of the Hochschild
(co)homology of the twisted group algebra, and
the third proposition says
that one can express the Hochschild homology of
the quantum algebra $A$ with values in the bimodule
$Ag$ twisted by the action of an element $g\in G$
in terms of De Rham cohomology of the subvarieties
of fixed points of $X$ under the action of $g$.
The third section of the paper is devoted to the
proof of proposition \ref{FLT} and the forth
most technical section is devoted to the proof
of proposition \ref{ONO}. We conclude our paper
by giving an application of our result to
deformation theory, namely to the definition of a
``global'' analog of symplectic reflection algebras
from \cite{EG}.

{\bf Notation.} Throughout the paper summation over repeated indices
is assumed. $X$ is a smooth affine algebraic variety (over $\bbC$)
equipped with an algebraic symplectic form.
The dimension of $X$ is $2n$\,. We omit symbol $\wedge$ referring to a local basis
of exterior forms as if we thought of $dx^i$'s as anti-commuting variables.
By a nilpotent linear operator we always mean
one whose second power vanishes. If an associative algebra
$A$ is endowed with an action of a finite group $G$ then by
$a^g$ we denote the left action of an element $g\in G$ on an
element $a\in A$ and by $\AG= \AGP$ we denote the corresponding
twisted group algebra. If $B$ is an associative algebra
then $B^{op}$ denotes the algebra $B$ with the opposite
multiplication. For a vector bundle $V$ we denote by
$S^m V$ the $m$-th symmetric tensor power of $V$.

{\bf Acknowledgments.} We would like to thank R. Anno, V. Ginzburg,
D. Kaledin, and A. Oblomkov for useful discussions.
This work is partially supported by the NSF grant
DMS-9988796, the Grant for Support of Scientific
Schools NSh-1999.2003.2 and the grant
CRDF RM1-2545-MO-03.

\section{Sketch of the proof.}
We start with an observation that $A$ is a simple algebra.
Indeed, if $I$ is a non-zero two-sided ideal in $A$ then
$I_+=I\cap A_+$ is a saturated ideal in $A_+$. This implies that
$I_0=I_+/\h I_+$ is a non-zero Poisson
ideal in $A_0={\cal O}(X)$. But $X$ is symplectic and
therefore $I_0=A_0$\,. Since $I_+$ is saturated, $I_+=A_+$
and hence $I=A$\,.

Next, we would like to mention the following
\begin{pred}[\cite{Mont}]
If $B$ is a simple algebra and
$G$ is a finite group which acts on $B$
by outer automorphisms then the twisted
group algebra $B[G]$ is also simple.
\end{pred}
It is not hard to see that the algebra $A$ does not have
inner automorphisms of finite order. Therefore,
since $G$ acts on $A$ faithfully, we conclude
that the twisted group algebra $\AG$ is simple.
The latter implies that the $(A^G\,, \AG)$-bimodule
$$
P_{\bf e}= {\bf e} \AG
$$
and the $(\AG\,, A^G)$ -bimodule
$$
Q_{\bf e} = \AG {\bf e}\,,
\qquad {\bf e} = \frac1{|G|}\sum_{g\in G} g
$$
establish a Morita equivalence between $\AG$ and $A^G$\,.
Hence, we can replace $HH^{\bul}(A^G)$ in (\ref{main})
by the graded vector space $HH^{\bul}(\AG)$ of Hochschild
cohomology of the twisted group
algebra $\AG$.

Let $g$ be an element of the group $G$.
The subvariety $X_g$ of fixed points in $X$ under the action
of $g$ may, in principle, be reducible and have components of
different dimensions. However,
it is not hard to show that $X_g$ is smooth\footnote{To prove this
assertion one can use normal coordinates associated with a $G$-invariant
connection on $X$.}. Therefore,
if $\{ X^Q_g \}\,,~Q=1,2, \dots $ is the set of irreducible components of $X_g$
then two different components do not intersect each other.
Thus each irreducible (and also connected) component $X^Q_g$ of
$X_g$ is a smooth affine subvariety of $X$.

For the orbifold modeled on the quotient $X/G$ of a smooth affine
symplectic variety $X$ (over $\bbC$) the definition of the Chen-Ruan cohomology
\cite{CR} simplifies drastically. Namely, we have \cite{GK}
that
\begin{equation}
\label{Ch-R}
H^{\bul}_{CR}(X/G) = \left(\bigoplus_{g\in G, Q}
H^{\bul-codim X^Q_g}_{DR}(X^Q_g) \right)^G\,,
\end{equation}
where the notation $H_{DR}$ stands for algebraic De Rham cohomology.
Notice that our coefficient field is $\bbC((\h))$. Hence,
by the algebraic De Rham cohomology of a variety $Y$
we mean the cohomology of the De Rham complex of the formal
Laurent power series in $\h$ with values in algebraic exterior
forms on $Y$. By Grothendieck's theorem this
cohomology is isomorphic to the singular cohomology of $Y$
with the coefficients in the field $\bbC((\h))$\,.

The proof of theorem \ref{ONA} is essentially
based on the three propositions given below.

\begin{pred}[M. Van Den Bergh,\cite{VB}]
\label{VdB} Suppose that an algebra $B$ over a field $k$ has
a finite Hochschild dimension and
$$
Ext^i_{B\otimes B^{op}}(B, B\otimes B^{op})=
\cases{
\begin{array}{lr}
U & {\rm if}~ i = d\,, \\
0 & {\rm otherwise}
\end{array}
}
$$
as $B$-bimodules, where $U$ is an invertible
$B$-module. Then for any $B$-bimodule $M$ the
graded vector spaces
$$
HH^{\bul}(B, M) \qquad {\rm and} \qquad HH_{d - \bul}(B,U\otimes_{B} M)
$$
are naturally isomorphic.
\end{pred}
Following \cite{EO} we denote by $VB(d)$ the
class of associative algebras $B$ satisfying
the conditions of proposition \ref{VdB} with
$U=B$.

Notice that if $M$ is a bimodule over the twisted group
algebra $\AG$ then $M$ is also a $G$-bimodule and an
$A$-bimodule. Therefore, to any $\AG$-bimodule $M$ we
assign the Hochschild complexes $C_{\bul}(A, M)$
and $C^{\bul}(A,M)$ endowed with the following
action of $G$:
\begin{equation}
\label{action-G}
\al(g) (m, a_1,a_2, \dots, a_q) =
(g m g^{-1}, (a_1)^g, (a_2)^g, \dots, (a_q)^g)\,,
\end{equation}
\begin{equation}
\label{action-G1}
\al(g) (\Psi) (a_1,a_2, \dots, a_q) =
g\Psi((a_1)^{g^{-1}}, (a_2)^{g^{-1}}, \dots, (a_q)^{g^{-1}})g^{-1}\,,
\end{equation}
$$
g \in G\,, ~(m, a_1, \dots, a_q)\in C_{q}(A,M)\,,
~ \Psi \in C^{q}(A,M)\,.
$$

\begin{pred}
\label{FLT}
If $A$ is an associative algebra (over $\bbC$) acted
on by a finite group $G$, and $M$ is a bimodule over
the twisted group algebra $\AG$, then
\begin{equation}
\label{Lorenz-hom}
HH_{\bul} (\AG,M) = \Big( HH_{\bul} (A,M) \Big)^G\,,
\end{equation}
and
\begin{equation}
\label{Lorenz-coh}
HH^{\bul} (\AG,M) = \Big( HH^{\bul} (A,M) \Big)^G\,,
\end{equation}
where the action of $G$ on
$HH_{\bul}(A, M)$ (resp. on $HH^{\bul}(A,M)$)
is induced by (\ref{action-G})
(resp. by (\ref{action-G1})).
\end{pred}
\begin{pred}
\label{ONO}
The graded $G$-modules
$$
HH_{\bul} (A, \AG)
$$
and
$$
\bigoplus_{g\in G,Q} H^{dim\, X^Q_g-\bul}_{DR}(X^Q_g)
$$
are isomorphic.
\end{pred}
{\bf Remark.} We would like to mention that in
the case of the trivial group $G=\{e\}$ the above
proposition reduces to the algebraic analogue of
the well-known theorem of R. Nest and B. Tsygan
(theorem A$2.1$ in \cite{NT}) which is proved
for the quantum algebra of compactly supported
functions of a smooth symplectic manifold.

Theorem \ref{ONA} will follow immediately from the
above propositions if we prove that $\AG$ belongs to the
class $VB(2n)$, where $2n$ is the dimension of $X$.

For this purpose we first prove that
\begin{pred}
If an associative algebra $A$ (over $\bbC$)
satisfies the conditions of proposition \ref{VdB}
and
\begin{equation}
\label{vazhno}
Ext^d_{A\otimes A^{op}}(A,A\otimes A^{op})=A
\end{equation}
both as an $A$-bimodule and a $G$-module,
then the twisted group algebra $\AG$ also
belongs to the class $VB(d)$.
\end{pred}
{\bf Proof.} By the assumption $A$ has a finite Hochschild
dimension. Hence, due to proposition \ref{FLT}, so does
the algebra $\AG$.

Using proposition \ref{FLT} once again we get
$$
Ext^i_{\AG\otimes \bAG} (\AG, \AG \otimes \bAG) =
HH^i(\AG,\AG \otimes \bAG) =
$$
$$
(HH^i(A, \AG \otimes \bAG))^G =
\Big ( \bigoplus_{g,h \in G}
HH^i(A,Ag\otimes h A^{op})
\Big)^G\,.
$$
For any pair $g,h \in G$ the $A$-bimodule $A g \otimes hA^{op}$
is isomorphic to $A\otimes A^{op}$ and hence
$$
Ext^i_{\AG\otimes \bAG} (\AG, \AG \otimes \bAG) =0
$$
for any $i \neq d$.
Furthermore, using condition (\ref{vazhno}) we
get
$$
Ext^d_{\AG\otimes \bAG} (\AG, \AG \otimes \bAG) =
\Big ( \bigoplus_{h \in G}
\AG \otimes h \Big)^G\,,
$$
where $\AG\otimes h$ denotes the $\AG$-bimodule
$\AG$ twisted from the right by $h$, and
$g\in G$ acts on the direct sum
$\bigoplus_{h \in G} \AG \otimes h$ multiplying
the first component by $g^{-1}$ from the right
and the second component by $g$ from the left.

The direct sum
$\bigoplus_{h \in G} \AG \otimes h$ is
free as a $G$-module. Hence,
$$
Ext^d_{\AG\otimes \bAG} (\AG, \AG \otimes \bAG) = \AG
$$
as an $\AG$-bimodule and the desired statement follows. $\Box$

Due to the above proposition it suffices to prove that
\begin{pred}[\cite{EO}]
\label{P+L}
The algebra $A$ of quantum functions on $X$
has a finite Hochschild dimension,
$$
Ext^{i}_{A\otimes A^{op}}(A,A\otimes A^{op})=0
$$
for any $i\neq 2n=dim X$ and
\begin{equation}
\label{need}
Ext^{2n}_{A\otimes A^{op}}(A,A\otimes A^{op})=A
\end{equation}
both as an $A$-bimodule and as a $G$-module.
\end{pred}
The proof of this proposition given in \cite{EO}
misses some details.
For this reason we decided to repeat the proof
of \cite{EO} below and spend more time on the details
omitted in \cite{EO}.

~\\
{\bf Proof.} Let us suppose for a moment that
\begin{equation}
\label{assume}
A_0=\cO(X)\in VB(2n)
\end{equation}
and
\begin{equation}
\label{G-uni}
Ext_{A_0\otimes A_0}^{2n}(A_0, A_0\otimes A_0)= A_0
\end{equation}
as $G$-modules.

These assumptions have a pure algebraic geometric
meaning and they have nothing to do with the quantization
of $X$. For this reason we postpone the proofs
of the assumptions until the very end of
the section.

Since $A$ is a quantization of the commutative
algebra $A_0$, for any $A$-module $M$ the Hochschild
complex $C_{\bul}(A,M)$ is endowed with the filtration
induced by degrees in $\h$. For the corresponding
spectral sequence
$$
E^1_{p,q}(C_{\bul}(A,M)) = HH_{p+q}(A_0,M).
$$
Hence due to our assumption (\ref{assume})
$A$ has a finite Hochschild dimension.

Let us now consider the Hochschild complex
$C^{\bul}(A,A\otimes A^{op})$ of $A$-bimodules and
$G$-modules. The natural filtration
induced by degrees in $\h$ gives us the spectral sequence
converging to $Ext^{i}_{A\otimes A^{op}}(A,A\otimes A^{op})$
and having
\begin{equation}
\label{o-kak!}
E^1_{p,q}(C^{\bul}(A,A\otimes A^{op})) =
HH^{p+q}(A_0, A_0\otimes A_0)\,.
\end{equation}
Due to our assumptions (\ref{assume}) and (\ref{G-uni})
the right hand side of (\ref{o-kak!}) is isomorphic to $A_0$ as an
$A_0$-bimodule and as a $G$-module if $p+q=2n$, and
vanishes if $p+q\neq 2n$. The latter implies that
all the differentials
$d_1,~d_2,~d_3, ~\dots$
starting from $d_1$ in $E^{1}_{\bul, \bul}$
vanish. Thus,
$$
E^1_{p,q}(C^{\bul}(A,A\otimes A^{op})) =
E^{\infty}_{p,q}(C^{\bul}(A,A\otimes A^{op}))\,,
$$
$$
Ext^{i}_{A\otimes A^{op}}(A,A\otimes A^{op})=0
$$
for any $i\neq 2n$ and the associated graded
$A_0$-bimodule and $G$-module of the filtered
$A$-bimodule and $G$-module
\begin{equation}
\label{U-A}
U_A=Ext^{2n}_{A\otimes A^{op}}(A,A\otimes A^{op})
\end{equation}
is naturally isomorphic to
$A_0((\h))$.
This implies that there exists
a $G$-invariant automorphism
$$
\nu= Id+ \h \nu_1 + \dots
$$
of $A_+$ such that $U_A$
is isomorphic to $A\nu$ both as an $A$-bimodule
and as a $G$-module.

In particular, we see that $A$ satisfies the
conditions of proposition \ref{FLT}. Hence,
$$
HH^0(A, \nu^{-1}A) = HH_{2n}(A)\,.
$$
Using proposition \ref{ONO} for the trivial
group $G={e}$ we get that $HH_{2n}(A)= H^0(X)((\h))$.
Hence, $HH^0(A, \nu^{-1}A)=\bbC((\h))$ and therefore
there is a non-zero element $b\in A\,,~b\neq 0$
such that for any $a\in A_+$
$$
\nu^{-1}(a) \star b = b \star a
$$
or equivalently,
\begin{equation}
\label{nu-1}
a \star b = b \star \nu(a)\,, \qquad \forall~a\in A_+\,.
\end{equation}
We can safely assume that $b\in A_+$.

Denoting by $a_0$ and $b_0$ the zero-th terms
of the elements $a$ and $b$, respectively,
and comparing the terms of the first degree
in $\h$ in (\ref{nu-1}), we get
\begin{equation}
\label{nu-11}
\nu_1(a_0)b_0 = \{a_0,b_0\}\,.
\end{equation}
If $b_0$ vanishes at some point $p\in X$
then one can pick a function $a_0\in A_0=\cO(X)$
such that the left hand side of (\ref{nu-1}) has a bigger order of
vanishing than the right hand side. From this, we conclude
that $b_0$ is a nowhere vanishing function and therefore $b$ is an
invertible element of $A_+$. Hence equation (\ref{nu-1}) implies that
the automorphism $\nu$ is inner and $A$ indeed
belongs to the class $VB(2n)$.

To prove that (\ref{need}) is an isomorphism
of $G$-modules it suffices to show that the element
$b$ in (\ref{nu-1}) is $G$-invariant.
Since the automorphism $\nu$ is $G$-invariant we have that
for any $a\in A_+$ and $g\in G$
$$
b^{-1} \star a \star b = (b^g)^{-1} \star a \star b^g
$$
or equivalently,
$$(b^g \star b^{-1})\star a = a\star (b^g \star b^{-1})\,.$$
But the algebra $A_+$ has a one-dimensional center
$\bbC[[\h]]\subset A_+$. Hence, $b$ generates a one-dimensional
representation of the group $G$ over the ring
$\bbC[[\h]]$:
$$
b^g =\al_g b\,, \qquad \al_g \in \bbC[[\h]]\,.
$$
Since the associated graded module of
the filtered $G$-module (\ref{U-A}) is
naturally isomorphic to
$A_0((\h))$
as a $G$-module, for any $g\in G$
$$
\al_g\Big|_{\h=0}=1\,.
$$
But we also have that $(\al_g)^{N_g}=1$, where
$N_g$ is the order of the element $g$. Hence,
$\al_g=1$ for any $g\in G$\,. $\Box$

We are now left with our assumptions (\ref{assume})
and (\ref{G-uni}) about the algebra $A_0=\cO(X)$
of regular functions on the affine variety $X$.
It is not hard to see that these assumptions are
consequences of the following proposition and the fact that
$X$ is endowed with a $G$-invariant algebraic
symplectic form.

\begin{pred}
\label{Ext-A-0}
If $A_0$ is the algebra of regular functions
on a smooth affine algebraic variety $X$ (over
a field of characteristic zero) of dimension $d$,
then $A_0$ has a finite Hochschild
dimension,
$$
Ext^i_{A_0\otimes A_0}(A_0, A_0\otimes A_0) = 0\,,
$$
if $i\neq d$, and the $A_0$-bimodule
\begin{equation}
\label{A0-VB}
Ext^d_{A_0\otimes A_0}(A_0, A_0\otimes A_0)
\end{equation}
is naturally isomorphic to
$$
Hom_{A_0}(\Om^d(A_0),A_0)\,,
$$
where $\Om^d(A_0)$ is the module of
top degree exterior forms on $X$.
\end{pred}
{\bf Proof} is parallel to the proof of the
Hochschild-Kostant-Rosenberg theorem \cite{HKR}\,.
It is based on the standard lemma of
commutative algebra which says that if
a morphism
$$
\mu \,:\,M_1 \mapsto M_2
$$
of modules over a commutative algebra
$B$ becomes an isomorphism after the localization
with respect to any maximal ideal $\mm$ of $B$
then $\mu$ is an isomorphism of modules.

First, since $X$ is smooth, the kernel
of the multiplication map
$$
\mu_0\,:\, A_0 \otimes A_0 \mapsto A_0
$$
is a locally complete intersection.
In other words, for any maximal ideal $\mm\subset A_0$,
the kernel of the localized map
$$
(\mu_0)_{\mm}\,:\, (A_0 \otimes A_0)_{\mu^{-1}(\mm)}
\mapsto (A_0)_{\mm}
$$
is generated by a regular sequence of length $d$\,.

Thus, applying to $(A_0)_{\mm}$ the Koszul resolution,
we get that for any maximal ideal $\mm\subset A_0$
and for any $i\neq d$
$$
(Ext^i_{A_0\otimes A_0}(A_0, A_0\otimes A_0))_{\mm} =0\,.
$$
Hence,
$$
Ext^i_{A_0\otimes A_0}(A_0, A_0\otimes A_0) =0
$$
for any $i\neq d$\,.

If $i=d$ then we have a pair of morphisms of
modules
$$
Ext^d_{A_0\otimes A_0}(A_0, A_0\otimes A_0) \mapsto
Ext^d_{A_0\otimes A_0}(A_0, A_0)\,,
$$
\begin{equation}
\label{pair}
Hom_{A_0}(\Om^d(A_0),A_0)
\mapsto
Ext^d_{A_0\otimes A_0}(A_0, A_0)\,.
\end{equation}
The first morphism is induced by multiplication
$\mu_0$ and the second one is the so-called
anti-symmetrization morphism given by the
formula
\begin{equation}
\label{HKR}
\ve(\ga)(a_1, \dots, a_d) =
\ga(d\,a_1\wedge \dots \wedge d\,a_d)\,,
\qquad
\ga\in Hom_{A_0}(\Om^d(A_0),A_0)\,.
\end{equation}
Using the Koszul resolution once again, one
can show that for any maximal ideal $\mm\subset A_0$
the localizations of the morphisms (\ref{pair}) are
isomorphisms. Thus the above lemma from commutative
algebra  completes the proof of the
proposition. $\Box$

\section{The decomposition of
the (co)homology of the twisted group algebras.}

{\bf Proof} {\it of proposition} \ref{FLT}.
The idea of the proof is to guess a free left
resolution of $\AG$ as an $\AG\otimes \AG^{op}$-module which allows
us to get the desired decompositions (\ref{Lorenz-hom})
and (\ref{Lorenz-coh}) for free. Namely, we consider
the following complex of free
$\AG\otimes \AG^{op}$-modules
$$
\begin{array}{ccccccccc}
 ~&
\downarrow^{\beta'}&
  ~&
  \downarrow^{\beta'} & ~&\dots&
   ~&  ~& ~\\[0.3cm]
 \stackrel{\beta}{\rightarrow}&
C_{m,q}&
  \stackrel{\beta}{\rightarrow}&
  C_{m-1,q}&
  \stackrel{\beta}{\rightarrow}& \dots&
   ~& ~ & ~ \\[0.3cm]
 ~&
\downarrow^{\beta'}&
  ~&
  \downarrow^{\beta'} & ~&\dots&
   ~&  ~& ~ \\[0.3cm]
 \stackrel{\beta}{\rightarrow}&
C_{m,q-1}&
  \stackrel{\beta}{\rightarrow}&
  C_{m-1,q-1}&
   \stackrel{\beta}{\rightarrow}& \dots&
   ~&  ~& ~\\[0.3cm]
 ~&
\downarrow^{\beta'}&
  ~&
  \downarrow^{\beta'} & ~&\dots&
   ~&  ~& ~\\[0.3cm]
~&
\dots&
  ~&
  \dots & \dots& ~&
\dots& ~& ~\\[0.3cm]
~&
~&
  ~&
  ~ & ~&~&
   ~&  ~& ~\\[0.3cm]
~&
~&
  ~&
  ~ &  \downarrow^{\beta'}&~&
    \downarrow^{\beta'}&  ~&  ~ \\[0.3cm]
 ~&~& \dots& \stackrel{\beta}{\rightarrow} &
   \AG\otimes A\otimes \AG\otimes k[G] & \stackrel{\beta}{\rightarrow} &
   \AG\otimes \AG\otimes k[G]& \stackrel{\beta}{\rightarrow} &
   0 \\[0.3cm]
 ~&~&~&~&
   \downarrow^{\beta'} & ~ &
     \downarrow^{\beta'}  & ~ &
      ~ \\[0.3cm]
 ~&~&\dots& \stackrel{\beta}{\rightarrow}&
   \AG\otimes A\otimes \AG& \stackrel{\beta}{\rightarrow} &
   \AG\otimes \AG& \stackrel{\beta}{\rightarrow} & \AG\,,
\end{array}
$$
where
\begin{equation}
\label{C-dot-dot}
C_{m,q}=\AG\otimes A^{\otimes m}\otimes \AG\otimes k[G]^{\otimes q}\,.
\end{equation}
The horizontal differential $\beta$ is defined by
$$
\beta (ag\otimes a_1\otimes \dots \otimes
a_m\otimes bh\otimes g_1\otimes\dots \otimes g_q)=
$$
\begin{equation}
\label{beta-a}
(-)^q (aa^g_1g \otimes a_2 \dots
\otimes a_m\otimes bh\otimes g_1\otimes\dots \otimes g_q -
ag \otimes a_1a_2 \otimes a_3 \otimes \dots
\otimes a_m\otimes bh\otimes g_1\otimes\dots \otimes g_q +
\end{equation}
$$
\dots +(-)^{m} ag\otimes a_1 \otimes \dots
\otimes a_m bh\otimes g_1\otimes\dots \otimes g_q)
$$
where $a, a_1, \dots, a_m, b\in A$ and
$g,h,g_1, \dots, g_q\in G$.
If $m=0$ and $q\neq 0$ then $\beta$ vanishes and
if both $m=0$ and $q=0$ then $\beta$ is the ordinary
multiplication in $\AG$\,.

The vertical differential $\beta'$ is just
the boundary operator of the group Eilenberg-Mac Lane
complex:
$$
\beta' (ag\otimes a_1\otimes \dots \otimes
a_m\otimes bh\otimes g_1\otimes\dots \otimes g_q)=
$$
$$
agg_1\otimes a^{g^{-1}_1}_1 \otimes a^{g^{-1}_1}_2 \dots
\otimes a^{g^{-1}_1}_m\otimes b^{g^{-1}_1} g^{-1}_1h\otimes
g_2\otimes\dots \otimes g_q
$$
\begin{equation}
\label{beta'-a}
- ag \otimes a_1\otimes \dots
\otimes a_m\otimes bh\otimes g_1g_2\otimes g_3\otimes\dots \otimes g_q +
\end{equation}
$$
\dots +(-)^{q-1} ag\otimes a_1 \otimes \dots
\otimes a_m\otimes bh\otimes g_1\otimes\dots g_{q-2}\otimes g_{q-1}
g_q
$$
$$
+(-)^q ag\otimes a_1 \otimes \dots
\otimes a_m\otimes bh\otimes g_1\otimes\dots g_{q-2}\otimes g_{q-1}
$$
where $a, a_1, \dots, a_m, b\in A$ and
$g,h,g_1, \dots, g_q\in G$.
It is not hard to see that both differentials respect
the natural  $\AG\otimes \AG^{op}$-module structure:
$$
a_0g_0(ag\otimes a_1\otimes \dots \otimes
a_m\otimes bh\otimes g_1\otimes\dots \otimes g_q)b_0 h_0=
$$
\begin{equation}
\label{B-B-mod}
a_0a^{g_0}g_0g\otimes a_1\otimes \dots \otimes
a_m\otimes bb^{h}_0hh_0\otimes g_1\otimes\dots \otimes g_q,
\end{equation}
and anti-commute with each other:
$$
\beta \beta' + \beta' \beta=0\,.
$$

It turns out that the total complex of
(\ref{C-dot-dot}) is acyclic in all terms. We prove this
by defining explicitly the homotopy operator $\chi$
which acts on chains as follows. If $m\neq 0$ then
\begin{equation}
\label{hom-op}
\chi(ag\otimes a_1\otimes \dots \otimes
a_m\otimes bh\otimes g_1\otimes\dots \otimes g_q)=
(-)^q g\otimes a^{g^{-1}}\otimes a_1\otimes \dots \otimes
a_m\otimes bh\otimes g_1\otimes\dots \otimes g_q\,,
\end{equation}
and if $m=0$ then
\begin{equation}
\label{hom-op0}
\chi(ag\otimes bh\otimes g_1\otimes\dots \otimes g_q)=
(-)^q g\otimes a^{g^{-1}}\otimes bh\otimes g_1\otimes\dots \otimes g_q
+1 \otimes ab^g gh \otimes g\otimes g_1 \dots g_q\,.
\end{equation}
Finally, the action of $\chi$ on $\AG$ is defined by
\begin{equation}
\label{hom-op00}
\chi (ag) = 1\otimes ag\,.
\end{equation}
The desired property of $\chi$
\begin{equation}
\label{hom-op-pro}
\chi (\beta+\beta') + (\beta+\beta')\chi = Id
\end{equation}
can be proved by the direct computation.

Since our base field $\bbC$ obviously contains
the inverse $|G|^{-1}$ of the order of the group $G$,
the rest is evident from the definition of Hochschild
(co)homology, the isomorphism between the module of invariants
and the module of coinvariants, and the acyclicity
of the (co)homological group Eilenberg-Mac Lane complex
in positive dimensions. $\Box$

{\bf Remark.} Different statements closely related to the
assertion of proposition \ref{FLT} were obtained by
various authors. Thus, the decomposition
(\ref{Lorenz-hom}) was given in preprint \cite{J-L} of J.-L. Brylinski.
In the same year J.-L. Brylinski \cite{J-L1} considered
cyclic homology of a natural generalization of
the twisted group  algebra to the case when the
initial algebra is the algebra of functions on a
smooth manifold and the manifold is acted by a Lie group.
J.-L. Brylinski suggested in \cite{J-L1} that the cyclic
homology of this algebra is a good replacement of
the equi\-va\-riant $K$-theory in the case when
the Lie group is noncompact.
Simultaneously with J.-L. Brylinski, B. Feigin and B. Tsygan
proposed in \cite{FT} the corresponding
decomposition for the cyclic homology of an
arbitrary twisted group algebra.  Another derivation
of this result of  B. Feigin and B. Tsygan
was obtained in paper \cite{GJ} by
E. Getzler and J.D.S. Jones. The proof given
in \cite{GJ} involves a notion of the
paracyclic module which is a natural generalization of
the cyclic module. Using this notion in
\cite{GJ1}, J. Block, E. Getzler and  J. D. S. Jones considered the
case of the topological algebra acted by a Lie group.
Finally, in \cite{Lorenz} M. Lorenz proved the
decomposition (\ref{Lorenz-hom}) of the Hochschild
homology for any $G$-graded associative algebra
which is a natural generalization of
the twisted group algebra.

\section{Hochschild homology $HH_{\bul}(A, Ag)$.}
\subsection{The case of the symplectic affine space.}
Before speaking about the general symplectic affine
variety, we will mention the necessary results for
the case of the symplectic affine space acted on
by the cyclic group $\{1,g,g^2, \dots, g^{N-1}\}$
of order $N$. These results we need are essentially contained
in paper \cite{Alev}. The only difference is that
unlike the authors of \cite{Alev},
we speak about the formal Weyl
algebra, rather than the ordinary one.

Let $V$ be a $2n$-dimensional vector space (over $\bbC$)
endowed with a symplectic form $B$.
Let $\{y^1, \dots, y^{2n}\}$ be a basis
in $V$. We denote by $||B^{ij}||$,
$(i,j=1, \dots , 2n)$ the matrix
$$
B^{ij} = B(y^i, y^j)
$$
of the form $B$ in this basis.

\begin{defi}
\label{fWa}
The Weyl algebra\footnote{More precisely, we should call it
{\it the formal Weyl algebra}.}  $W$ associated with the symplectic
vector space $V$ is the vector space $\bbC[[V]]((\h))$
of the formally completed symmetric algebra of $V$
equipped with the following
(associative) multiplication
\begin{equation}
\label{circ}
(a \circ b)(y,\h)=
\exp\left(\frac{\h}{2}B^{ij}\frac{\pa}{\pa y^i}\frac{\pa}{\pa z^j}
\right)a(y,\h)b(z,\h)\Big|_{y=z}.
\end{equation}
\end{defi}
One can easily see that the multiplication
defined by (\ref{circ}) does not depend on the
choice of a basis in $V$. We view $W$ as an algebra
over the field $\bbC((\h))$\,.

The Weyl algebra $W$ is naturally filtered with respect
to the degree of monomials $2[\h]+[y]$ where
$[\h]$ is the degree in $\h$ and $[y]$ is the degree in $y$:
\begin{equation}
\label{filtr}
\begin{array}{c}
\displaystyle
\dots \subset W^1  \subset W^0  \subset W^{-1}
\dots \subset W\,, \\[0.3cm]
\displaystyle
W^m = \{ a = \sum_{2k+p\ge m}
\h^k a_{k; i_1 \dots i_p} y^{i_1} \dots y^{i_p}\}
\end{array}
\end{equation}
This filtration defines the $2[\h]+[y]$-adic topology in $W$.

Let us assume that a symplectic vector space $(V,B)$ is
acted on by a cyclic group $Z_{N}=\{1,g,g^2, \dots, g^{N-1}\}$
of order $N$, and $g$ preserves the symplectic form
$B$. It is not hard to show that in this situation
$V$ becomes a direct sum of two symplectic
vector spaces:
\begin{equation}
\label{split}
V=  Ker\,(\pi(g)-Id) \oplus Im \,(\pi(g)-Id)\,,
\end{equation}
where $\pi(g)$ is the action of $g$ on $V$\,.
We denote the subspace $Ker\, (\pi(g)-Id)$ by $V_g$,
the dimension of $V_g$ by $2m_g$, and
the dual form to the corresponding
symplectic form on $V_g$ by $B^g$
$$
B^g = \frac1{2}B^g_{ab} y^a\wedge y^b\,,
$$
where $\{y^1, \dots, y^{2m_g}\}$ is a basis
in $V_g$\,.

The group $Z_N$ naturally acts on the Weyl algebra $W$.
In particular, we can define a $W$-bimodule $Wg$
twisted from the right by the action of $g$\,.
For this bimodule we have the following
\begin{pred} With the above notations,
\label{ochen-nado}
\begin{equation}
HH_{q}(W, Wg) = \cases{
\begin{array}{lr}
\bbC((\h))\,, & {\rm if}~ q = 2m_g\,, \\[0.3cm]
0\,,  &  {\rm otherwise}\,.
\end{array}
}
\end{equation}
If $\ve_{a_1\dots a_{2m_g}}$ are components of
the Liouville volume form
\begin{equation}
\label{ve}
\ve= \wedge^{m_g} B^g
\end{equation}
in the basis $\{ y^1, \dots, y^{2m_g}\}$
then the element
\begin{equation}
\label{ON}
\psi = \ve_{a_1 \dots a_{2m_g}} 1\otimes y^{a_1} \otimes \dots \otimes y^{a_{2m_g}}
\end{equation}
is a non-trivial cycle of the normalized
Hochschild complex $\bC(W,Wg)$\,.
\end{pred}
{\bf Proof.} The first claim was essentially proved in
\cite{Alev}. The only subtlety we have to mention
is that unlike the authors of \cite{Alev}
we consider the formal rather than
the ordinary Weyl algebra. However, it is
not hard to verify that the arguments of \cite{Alev}
can be extended to the formal Weyl algebra as
well\footnote{See, for example, appendix B  of paper \cite{Hoch-DR}
in which this issue is examined for the cohomological
complex in the case of the trivial group action.}.

Let us mention that formula (\ref{ON})
does not depend on the choice of the
local basis in $V_g$ and therefore
the chain $\psi\in \bC_{2m_g}(W,Wg)$ is well-defined.
A simple computation shows that (\ref{ON})
is indeed a cycle in the normalized Hochschild
complex $\bC_{\bul}(W,Wg)$.

To prove that the cycle (\ref{ON}) is non-trivial we
consider the spectral sequence associated with the
filtration induced by degrees in $\h$. It is obvious that
$$
E^1_{p,q}(W,Wg) = HH_{p+q}(\bbC[[V]],\bbC[[V]]g)\,,
$$
where $\bbC[[V]]$ is viewed as a commutative algebra.

The splitting (\ref{split}) gives us the natural
projection $pr$ from the ring $\bbC[[V]]$
to the ring $\bbC[[V_g]]$. Using this projection
we define the following analogue of the
antisymmetrization map
\begin{equation}
\label{anti-sym}
\mu (a_0, a_1, \dots, a_q) = pr (a_0)
 d(pr(a_1)) \wedge \dots \wedge d (pr(a_q))\,,
\end{equation}
from the normalized Hochschild complex
$\bC_{\bul}(\bbC[[V]], \bbC[[V]]g)\otimes \bbC((\h))$
to the graded vector space
$$
\bbC[[V_g]]((\h))\otimes\, \wedge^{\bul}(V_g)\,.
$$
With the help of the standard Koszul resolution
of $\bbC[[V]]$ one can show that (\ref{anti-sym})
is a quasi-isomorphism from
$\bC_{\bul}(\bbC[[V]], \bbC[[V]]g)\otimes \bbC((\h))$
to the complex
\begin{equation}
\label{forms-Vg}
(\bbC[[V_g]]((\h))\otimes\, \wedge^{\bul}(V_g), 0)
\end{equation}
with the vanishing differential.

Thus,
$$
E^1_{p, q} = \bbC[[V_g]]\otimes\, \wedge^{p+q}(V_g)\,.
$$

It is not hard to show that
the differential
$$
d_1\,:\, E^1_{p,q} \mapsto E^1_{p+1,q-2}
$$
is the De Rham differential in (\ref{forms-Vg}) twisted by
the symplectic Hodge operator \cite{Br}.
Hence,
\begin{equation}
\label{E-1}
E^2_{p,q} =\cases{
\begin{array}{lr}
\bbC\,, & {\rm if}~ p+q = 2m_g\,, \\[0.3cm]
0\,,  &  {\rm otherwise}\,.
\end{array}
}
\end{equation}
It is easy to see that the cycle $\psi$ (\ref{ON})
projects onto a generator of $E^2$.
Therefore, $\psi$ is nontrivial and the proposition
follows. $\Box$

\subsection{Fedosov's construction.}
In this subsection we briefly recall Fedosov's
construction \cite{F} of the star-product on the symplectic
manifold $X$ focusing our attention on
the equivariant aspects of the construction.

The {\it Weyl algebra bundle} $\W$ is a bundle over $X$
whose sections are the following formal power series
\begin{equation}
a=a(x,y,h)=\sum_{k,l}\h^k a_{k; i_1 i_2\dots i_l}(x)y^{i_1}\dots
y^{i_l}\,,
\end{equation}
where $y=(y^1\dots y^{2n})$ are fiber coordinates of the tangent
bundle $TX$, $a_{k; i_1 i_2 \dots i_l}(x)$ are symmetric
covariant tensors, and the summation over $k$
is bounded below.

Multiplication of two sections of $\W$ is given by
the Weyl formula
\begin{equation}
\label{circ1}
a \circ b(x,y,\h)=\exp\left(\frac{\h}{2}\,
\om^{ij}\frac{\pa}{\pa y^i}\frac{\pa}{\pa z^j}\right)a(x,y,\h)b(x,z,\h)|_{y=z},
\end{equation}
where $\om^{ij}(x)$ are components of the corresponding Poisson
tensor
$$
\bom = \om^{ij}(x)\frac{\pa}{\pa x^i} \wedge \frac{\pa}{\pa x^j}\,,
\qquad
\om^{ik}\om_{kj} = \de^i_j\,.
$$

For any point $p\in X$ the fiber $\W_p$ of the Weyl algebra bundle
at $p$ is isomorphic to the Weyl algebra $W$ associated
with the cotangent space $T^*_p(X)$ at the point $p$ with the
symplectic form $\bom_p$\,.

The filtration (\ref{filtr}) of the Weyl algebra
gives us a natural filtration of the bundle $\W$
\begin{equation}
\label{filtr-bun}
\begin{array}{c}
\displaystyle
\dots \subset \W^1  \subset \W^0  \subset
\W^{-1} \dots \subset \W\,, \\[0.3cm]
\displaystyle
\G(\W^m) = \{ a = \sum_{2k+p\ge m}
\h^k a_{k;i_1 \dots i_p}(x)y^{i_1} \dots y^{i_p}\}\,.
\end{array}
\end{equation}
This filtration defines a $2[\h]+[y]$-adic topology in the
algebra $\G(\W)$ of sections of $\W$\,.

The vector space $\Om^{\bul}(\W)$ of smooth exterior
forms with values in $\W$ is naturally a graded
associative algebra with the product induced by (\ref{circ1}) and the
following graded commutator
$$
[a,b]= a \circ b - (-)^{q_a q_b} b\circ a\,,
$$
where $q_a$ and $q_b$ are exterior degrees of
$a$ and $b$\,, respectively.
The filtration of $\W$ (\ref{filtr-bun}) gives us a filtration
of the algebra $\Om^{\bul}(\W)$
$$
\begin{array}{c}
\dots \subset \Om^{\bul}(\W^1)
\subset\Om^{\bul}(\W^0)\subset  \phantom{aaaaaaa}  \\[0.3cm]
 \phantom{aaaaaaa} \subset\Om^{\bul}(\W^{-1}) \subset
\dots \subset \Om^{\bul}(\W)\,.
\end{array}
$$
(In the definition of this
filtration, the degree of $dx^i$ is put to be zero).
Notice that the algebra $\Om(\W)$ is endowed with
the natural action of the group $G$\,.

It is not hard to show that on $TX$ there exists
a torsion free connection $\pa^s$ compatible with
the symplectic structure $\om$ and invariant under
the action of $G$. Using this connection we
define the following linear operator
$$
\n \,:\,\Om^{\bul}(\W) \mapsto
\Om^{\bul+1}(\W)\,,
$$
\begin{equation}
\label{nabla}
\n=dx^i\frac{\pa}{\pa x^i}-
dx^i\G^{k}_{ij}(x)y^j\frac {\pa}{\pa y^k}\,,
\end{equation}
where $\G^k_{ij}(x)$ are the Christoffel symbols
of $\pa^s$\,.

Thanks to the compatibility of $\pa^s$ with the
symplectic structure $\om$ and with the action of
$G$ the operator (\ref{nabla})
is a $G$-equivariant derivation of the graded algebra
$\Om^{\bul}(\W)$. A simple
computation shows that
$$
\n^2 a = \frac1{2}[R,a]\,, \qquad
\forall~ a\in \Om(\W)\,,
$$
where
$$
R= \frac1{2\h} \om_{km}(R_{ij})^m_l(x) y^k y^l dx^i dx^j,
$$
and $(R_{ij})^m_l(x)$ is the Riemann curvature
tensor of $\pa^s$\,.

\begin{defi} The~ Fedosov~ connection~ is~ a~ nilpotent
derivation of the graded algebra $\Om^{\bul}(\W)$
of the following form
\begin{equation}
\label{F-conn}
D= \n + \frac1{\h}[\cA, \cdot\,]\,,
\qquad
\cA= - dx^i \om_{ij}(x)y^j + r\,,
\end{equation}
where $r$ is an element in
$\Om^1(\W^3)$
\end{defi}
The flatness of $D$ is equivalent to the fact
that the Fedosov-Weyl curvature
\begin{equation}
\label{Weyl-c}
C^W= \h R + 2 \n \cA + \frac1{\h} [\cA, \cA]
\end{equation}
of $D$ belongs to the subspace
$\Om^2(X)((\h))\subset \Om^2(\W)$\,.
A simple analysis of degrees in $\h$ and $y$
shows that $C^W$ is of the form
\begin{equation}
\label{nado}
C^W = -\om + \Om_{\h}\,,
\qquad \Om_h \in \h \Om^2(X)[[\h]]
\end{equation}
whereas the Bianchi identity $D (C^W)=0$ implies that
$\Om_h$ is a series of two-forms closed with
respect to the De Rham differential.

Let us remark that the Fedosov connection (\ref{F-conn})
can be rewritten as
\begin{equation}
\label{DDD}
D = \n -\de + \frac1{\h}[r, \cdot\,]\,,
\end{equation}
where
\begin{equation}
\label{de}
\de= \frac1{\h}[dx^i\om_{ij}(x)y^j, \cdot\,]=
dx^i \frac{\pa}{\pa y^i}
\end{equation}
is the Koszul derivation of the algebra
$\Om^{\bul}(\W)$\,.

For our purposes we will need
the homotopy operator for the Koszul
differential $\de$
\begin{equation}
\delta^{-1}a=y^k i\left( \frac \partial {\partial x^k}\right)
\int\limits_0^1 a(x,\h, ty,tdx)\frac{dt}t,  \label{del-1}
\end{equation}
where $i(\partial /\partial x^k)$ denotes the contraction
of an exterior form with the vector field
$ \partial /\partial x^k$\,, and $\delta ^{-1}$ is extended to
$\G(\W)$ by zero.

Simple calculations show that
$\de^{-1}$ is indeed the homotopy operator for $\de$,
namely
\begin{equation}
a=\sigma(a) +\delta \delta ^{-1}a + \delta ^{-1}\delta a\,,
\qquad \forall~a\in\Om(\W)
 \label{Hodge}
\end{equation}
where $\si$ is the natural projection
\begin{equation}
\label{sigma}
\si (a)= a \Big |_{y=0,~dx=0}\,, \qquad
a\in \Om^{\bul}(\W)
\end{equation}
from $\Om^{\bul}(\W)$ onto
the vector space $A_0((\h))$, where
$A_0=\cO(X)$\,.

The proof of the following theorem is contained in
section $5.3$ of \cite{F1}. More precisely, see theorem
$5.3.3$ and remarks at the end of section $5.3$.
\begin{teo}[Fedosov, \cite{F1}]
\label{F}
If $\pa^s$ is a symplectic connection on $X$ and $\Om_{\h}$
is a series of closed two-forms in $\h\Om^2(M)[[\h]]$ then
\begin{enumerate}
\item Iterating the equation
\begin{equation}
\label{iter}
r=\delta ^{-1}(R-\Om_{\h}) + \de^{-1}(\nabla r+
\frac 1{\h }r\circ r)
\end{equation}
in degrees in $y$'s
we get an element $r\in \Om^1(\W^3)$
such that the derivation
\begin{equation}
\label{DDD1}
D=\n - \de +\frac1{\h}[r, \cdot\,]
\end{equation}
is nilpotent and has the Fedosov-Weyl curvature (\ref{Weyl-c})
$C^W= -\om + \Om_{\h}$\,.

\item Given a Fedosov connection (\ref{F-conn})
one can construct a vector space isomorphism $\la$
\begin{equation}
\label{la}
\la\,:\, A_0((\h))~ \widetilde{\rightarrow}~ \G_D(\W)
\end{equation}
from $A_0((\h))$
to the algebra $\G_D(\W)$
of flat sections of $\W$ with respect
to $D$. The inverse isomorphism to $\la$
is the projection $\si$ (\ref{sigma})
and the product in $A_0((\h))$
\begin{equation}
\label{star}
a * b = \si(\la (a)\circ \la(b))\,, \qquad
a,b\in A_0((\h))
\end{equation}
induced via the isomorphism $\la$ is the
desired star-product in $A$.
\end{enumerate}
\end{teo}
{\bf Remark 1.} The cohomology class of the
Fedosov-Weyl curvature (\ref{Weyl-c}) is a well-defined
characteristic class of a star-product in $A_0((\h))$.
This characteristic class is referred to as
{\it the Deligne-Fedosov class}.\\
{\bf Remark 2.} Due to the naturality of the Fedosov construction
the Fedosov differential (\ref{DDD1}) corresponding to
a $G$-invariant series $\Om_{\h}$ of closed $2$-forms
is $G$-equivariant and the corresponding star-product
(\ref{star}) is $G$-invariant.\\
{\bf Remark 3.} For any function $a \in A_0=\cO(X)$
\begin{equation}
\label{la-der}
\left(\frac{\pa}{\pa y^{i_1}} \dots \frac{\pa}{\pa y^{i_k}}
\la(a)\right)\Big|_{y=0} =
 \pa_{x^{i_1}} \dots \pa_{x^{i_p}} a(x)+
{\rm lower~order~ derivatives~ of}~a\,.
\end{equation}

Due to the algebraic geometric version \cite{Hoch-DR} of
the theorem of P. Xu \cite{Xu} the star-product $\star$ in $A$ is
equivalent to some Fedosov star-product (\ref{star})
corresponding to the Deligne-Fedosov class
$[-\om + \Om_{\h}]$\,. We claim that
\begin{pred}
\label{star-ast}
The star-product $\star$ in $A$ is equivalent to
the Fedosov star-product (\ref{star}) corresponding to
the formal series $\Om_{\h}$ of $G$-invariant closed
$2$-forms and the equivalence
between $\star$ and (\ref{star}) can be established
by a $G$-invariant operator.
\end{pred}
{\bf Proof.} Due to theorem $2$ in \cite{Hoch-DR} the
$\star$ in $A$ is equivalent to the Fedosov star-product
corresponding to the formal series $\Om_{\h}$ such
that for any $g\in G$
$$
\pi(g)^*(\Om_h) - \Om_{\h}
$$
is a series of exact $2$-forms.
Therefore, the form
$$
\frac1{|G|}\sum_{g\in G}\pi(g)^{*}(\Om_{\h})
$$
is $G$-invariant, represents the same
class as $\Om_h$, and gives the Fedosov star-product
(\ref{star}) which is equivalent to $\star$ in $A$.

The second claim is proved in \cite{Hoch-DR}
(see corollary $1$). $\Box$

\subsection{The proof of proposition \ref{ONO}.}
The idea of the proof is based on the construction of
a complex of $G$-modules $\cK^{\bul}$ and a pair of
$G$-equivariant quasi-isomorphisms
\begin{equation}
\label{dir-sum}
\bigoplus_{g\in G} C_{\bul}(A,Ag)
\,\stackrel{?}{\longrightarrow}\,
\cK^{\bul}
\,\stackrel{?}{\longleftarrow}\,
\bigoplus_{g,Q} (\Om^{2m^Q_g - \bul}(X^Q_g)((\h)), d)\,,
\end{equation}
where $C_{\bul}(A,Ag)$ is the Hochschild complex
of $A$ with values in the twisted bimodule $Ag$,
$X^Q_g$ as above denotes the $Q$-th connected
(and irreducible) component of the subvariety
$X_g$ of the points fixed under the action
of $g \in G$, and $2m^Q_g=dim X^Q_g$\,.

Due to proposition \ref{star-ast} and remark $2$ after
theorem \ref{F} we may safely assume that the product
in $A$ is the Fedosov star-product (\ref{star}) corresponding
to a $G$-invariant series $\Om_{\h}$ of closed $2$-forms
and the Fedosov differential (\ref{DDD1}) is $G$-equivariant.

Let us denote by $E^{Q,g}$
the restriction of the tangent bundle
$TX$ to $X^Q_g$
$$
TX\Big|_{X^Q_g} = E^{Q,g}\,.
$$
Similarly we restrict the Weyl algebra bundle $\cW$ over
$X$ to $X^Q_g$ and denote this new bundle of algebras
by $\cW^{Q,g}$. The sections of $\cW^{Q,g}$ are the
following formal power series in fiber
coordinates $y^i$ of the bundle $E^{Q,g}$
\begin{equation}
\label{WQg-chain}
a(x,y, \h)=\sum_{k,p}
\h^k a_{k;i_1 \dots i_p}(x)y^{i_1} \dots y^{i_p}\,,
\end{equation}
where $a_{k;i_1 \dots i_p}(x)$ are components of
sections of $S^p(E^{Q,g})^*$ and the summation in $\h$
is bounded below.

Since any element $h\in G$ maps the subvariety $X^Q_g$
onto the subvariety $X^{Q'}_{hgh^{-1}}$, the sets of bundles
$\{E^{Q,g} \mapsto X^Q_g \}$ and $\{\cW^{Q,g} \mapsto X^Q_g \}$
are equipped with an action of $G$. In particular, the bundle $\cW^{Q,g}$
is acted on by the cyclic group $\ZNg$ generated by $g\in G$\,.

To the $\ZNg$-bundle $\cW^{Q,g}$ of algebras
we attach the bundle of fiberwise Hochschild chains
$\cC^{Q,g}$ over $X^Q_g$
\begin{equation}
\label{cC}
\cC^{Q,g} = \bigoplus_{m=0}^{\infty}\cC^{Q,g}_m\,,
\end{equation}
where
$$
\cC^{Q,g}_m = (\cW^{Q,g})^{\totimes~(m+1)}\,,
$$
and $\totimes$ stands for the tensor product
over $\cO(X^Q_g)((\h))$ completed in the
adic topology (\ref{filtr-bun}).

The sections of $\cC^{Q,g}_m$
are the following formal power series
in $m+1$ collections of fiber
coordinates $y^i_0, \dots, y_m^i$ of
$E^{Q,g}$
\begin{equation}
\label{f-chain}
a(x,\h, y_0, \dots, y_m)=\sum_k\sum_{\al_0 \dots \al_m}
\h^k a_{k;\al_0\dots \al_m}(x)y_0^{\al_0} \dots y^{\al_m}_m\,,
\end{equation}
where the summation in $k$ is bounded below
$\al$'s are multi-indices $\al={j_1\dots j_l}$,
$$
y^{\al}= y^{j_1}y^{j_2} \dots y^{j_l}\,,
$$
$a_{\al_0\dots \al_m}(x)$ are sections of
the bundle $S^{|\al_0|}E^{Q,g} \otimes \dots \otimes S^{|\al_m|}E^{Q,g}$,
and $|\al|$ is the length of the multi-index $\al$.

The differential in $\G(\cC^{Q,g})$ is given by the formula
$$
(\mb a)(x,\h, y_0, \dots, y_{m-1}) =
$$
\begin{equation}
\label{cC-mb}
S_{0}a(x,\h, z, \pi(g)(y_0), \dots, y_{m-1})\Big|_{z=y_0} -
S_{1}a(x,\h, y_0, z, y_1, \dots, y_{m-1})\Big|_{z=y_1} + \dots
\end{equation}
$$
+ (-)^{m-1} S_{m-1}a(x,\h, y_0, y_1, \dots, z, y_{m-1})\Big|_{z=y_{m-1}}
+(-)^m  S_0(x,\h, y_0, y_1, \dots, y_{m-1}, z)\Big|_{z=y_0}
$$
where $S_{a}$ ($a=0, \dots, m-1$) denotes the following operator:
$$
S_{a}= \exp
\left(\frac{\h}2\om^{ij}(x)\frac{\pa}{\pa z^i}\frac{\pa}{\pa y_a^i}\right)\,,
$$
and $\pi(g)$ stands for the action of $g$ on
the fibers of $E^{Q,g}$.

Notice that, since $\cW^{Q,g}$ is
the restriction of the Weyl algebra bundle $\cW$ to
$X^Q_g$, we can restrict the flat
Fedosov connection (\ref{DDD1}) on $\cW$ to
the flat connection on $\cW^{Q,g}$:
\begin{equation}
\label{F-conn-Qg}
D^{Q,g}= \n^{Q,g} + \frac1{\h}[\cA^{Q,g}, \cdot\,]\,,
\end{equation}
where $\n^{Q,g}$ and $\cA^{Q,g}$ are
the restrictions of the
connection (\ref{nabla}) and
form $\cA$ in (\ref{F-conn})
to $X^Q_g$, respectively.

Extending the connection $D^{Q,g}$ to the bundle
$\cC^{Q,g}_m = \cW^{\totimes~(m+1)}$ we get
the nilpotent differential on the graded vector
space $\Om^{\bul}(\cC^{Q,g}_{\bul})$ of exterior
forms on $X^Q_g$ with values in the sheaf
$\cC^{Q,g}_{\bul}$
\begin{equation}
\label{D-Qg}
D^{Q,g}\,:\,\Om^{\bul}(\cC^{Q,g}) \mapsto
\Om^{\bul+1}(\cC^{Q,g})\,.
\end{equation}
For the sake of simplicity, we use the same
notation $D^{Q,g}$ for this
differential.

Since the Fedosov differential $D$ (\ref{DDD1}) is
a $G$-equivariant derivation of the product (\ref{circ}),
the differential $D^{Q,g}$ anticommutes
with the fiberwise Hochschild differential $\mb$
$$
D^{Q,g} \mb + \mb D^{Q,g} =0
$$
if we use the convention
$$
\mb dx^i = -dx^i \mb\,.
$$

For our purposes we will also need the bundle $\bcC^{Q,g}$
of normalized fiberwise Hochschild chains over
each subvariety $X^Q_g$. Namely,
\begin{equation}
\label{bcC}
\bcC^{Q,g} = \bigoplus_{m=0}^{\infty}\bcC^{Q,g}_m\,,
\end{equation}
where
$$
\bcC^{Q,g}_m = \cW^{Q,g}\totimes(\overline{\cW}^{Q,g})^{\totimes~m}\,,
$$
$\totimes$ stands for the tensor product
over $\cO(X^Q_g)((\h))$ completed in the
adic topology (\ref{filtr-bun}), and
$\overline{\cW}^{Q,g}= \cW^{Q,g}/ \bbC((\h)) 1$\,.

It is not hard to see that the differentials
(\ref{cC-mb}) and (\ref{D-Qg}) still make sense for
the normalized chains (\ref{bcC}). Furthermore,
since the Fedosov differential (\ref{DDD1}) is
$G$-equivariant, the
direct sum of double complexes
\begin{equation}
\label{dir-s}
\cD_{p,q} = \bigoplus_{Q,g}
(\Om^{p}(\bcC_{q}^{Q,g}), D^{Q,g}+\mb)
\end{equation}
form a double complex of $G$-modules.

A standard lemma of homological algebra
(see, for example, proposition $1.6.5$ in
\cite{Loday}) says that the natural
projection
\begin{equation}
\label{normal}
\Pi\,:\,\bigoplus_{Q,g}
(\Om^{p}(\cC_{q}^{Q,g}), D^{Q,g}+\mb)
\mapsto \cD_{p,q}
\end{equation}
is a quasi-isomorphism of complexes.
We also mention that the projection (\ref{normal})
is $G$-equivariant.

We want to show that the total complex
$$
\cK^{\bul} = \bigoplus_{q-p=\bul } \cD_{p,q}
$$
of the double complex $\cD_{p,q}$ is
the desired middle term in (\ref{dir-sum}).

For this purpose we mention that the complex
$\bigoplus_{g\in G} C_{\bul}(A,Ag)$ of $G$-modules
can be naturally embedded into a bigger graded
$G$-module $J$
\begin{equation}
\label{Jets}
J= \bigoplus_{m} J_m, \qquad
J_m= \bigoplus_{g\in G, Q}\G(X^Q_g, Jets_{X^Q_g}(X^{m}))
\end{equation}
where $Jets_{X^Q_g}(X^{m})$ denotes the sheaf of jets
of regular functions on the product $X^m$ of $m$ copies
of $X$ near the subvariety $X^Q_g$ embedded into the
main diagonal $\D X^m \cong X$ of $X^m$.
It is not hard to see the Hochschild differential
of $C_{\bul}(A,Ag)$ still makes sense for the
graded module (\ref{Jets}) and is compatible
the action of $G$.

Although the complex (\ref{Jets}) is bigger
than the complex $\bigoplus_{g\in G} C_{\bul}(A,Ag)$,
we have
\begin{pred}
\label{to-formal}
The natural embedding
\begin{equation}
\label{q-iso}
\bigoplus_{g\in G} C_{\bul}(A, A g) \hookrightarrow J_{\bul}
\end{equation}
is a quasi-isomorphism of complexes.
\end{pred}
{\bf Proof.} It is clear from the construction that
(\ref{q-iso}) is a morphism of complexes of $G$-modules.
Furthermore, the map (\ref{q-iso}) induces a quasi-isomorphism
of $E_1$-terms of the spectral sequences associated
with the filtration induced by degrees in $\h$. Hence, due to
the standard snake-lemma argument of homological
algebra, the map (\ref{q-iso}) induces an isomorphism of
homology as graded $G$-modules. $\Box$

Notice that, thanks to remark $3$ after theorem \ref{F},
the map $\la$ (\ref{la}) gives us the natural embedding
of complexes of $G$-modules
\begin{equation}
\label{la-J}
\la^J\,:\, J_{\bul} \mapsto \bigoplus_{Q,g}\Om^0(\cC_{\bul}^{Q,g})\,.
\end{equation}
We claim that
\begin{pred}
\label{acyc-D}
For any $q$ the complex
$(\Om^{\bul}(\cC_{q}^{Q,g}), D^{Q,g})$ is
acyclic in positive dimension and the
morphism (\ref{la-J}) restricted to
$\G(X^Q_g, Jets_{X^Q_g}(X^{q}))$
gives us the isomorphism
\begin{equation}
\label{H0-D}
\la^{Q,g}\,:\,
\G(X^Q_g, Jets_{X^Q_g}(X^{q}))
\erarrow
Ker D^{Q,g} \cap \G (\cC^{Q,g}_q)\,.
\end{equation}
\end{pred}
{\bf Proof.} To prove the first claim we introduce,
on the complex $(\Om^{\bul}(\cC_{q}^{Q,g}), D^{Q,g})$,
the filtration induced by degrees in $y^a_0$ of
the sections (\ref{f-chain}), where $y^a_0$ are
fiber coordinates of the tangent bundle
$TX^Q_g\subset E^{Q,g}$. The zero-th
differential $d_0$ of the corresponding spectral
sequence has the form of the Koszul differential
\begin{equation}
\label{de-Qg}
\de^{Q,g} = dx^a \frac{\pa}{\pa y_0^a}\,,
\end{equation}
where $x^a$ are local coordinates on $X^Q_g$\,.
Hence, the $E^0$-term of the spectral sequence
is acyclic in positive dimension and the
first claim follows.

To prove the second claim we mention that
the morphism (\ref{H0-D}) is obviously
injective. Thus we are left with a proof
of surjectivity. In fact it suffices to
prove the surjectivity locally on $X^Q_g$\,.

Indeed, if $a$ is a $D^{Q,g}$-flat section
of $\cC^{Q,g}_q$ and for each local chart
$V\subset X^Q_g$ we have
$j_V\in \G(V, Jets_{X^Q_g}(X^{q}))$ such that
$a\Big |_{V}=\la^{Q,g} (j_V)$.
Then, by the injectivity of $\la^{Q,g}$, $j_V=j_V'$
on any intersection $V\cap V'$\,.
Hence we have $j\in Jets_{X^Q_g}(X^{q})$
such that $\la^{Q,g}(j)=a$\,.

An obvious analogue of theorem $5.5.1$ in \cite{F1}
says that the differential $D^{Q,g}$ on $\cW^{Q,g}$
can be locally conjugated to
$$
D^{Q,g}_0= dx^a \frac{\pa}{\pa x^a} -
dx^a \frac{\pa}{\pa y^a}\,,
$$
where $x^a$ are local coordinates on $X^Q_g$
and $y^a$ are fiber coordinates of the tangent bundle
$TX^Q_g\subset E^{Q,g}$\,. Thus the question of
surjectivity of $\la^{Q,g}$ reduces to a simple exercise of
the theory of partial differential equations. $\Box$

Next, we claim that
\begin{pred}
\label{el-kappa}
For any connected component $X^Q_g$
one can construct an element
$$
\ka^{Q,g}\in \Om^{0}(\bcC^{Q,g}_{2m^Q_g})\,,
$$
such that
\begin{equation}
\label{kappa}
D^{Q,g}\ka^{Q,g} = 0\,, \qquad
\mb \ka^{Q,g} =0\,,
\end{equation}
the evaluation $(\ka^{Q,g})\Big|_{p}$ at
any point $p\in X^Q_g$ is a non-trivial
$\mb$-cycle: in the fiber $(\bcC^{Q,g})_{p}$\,,
and for any $h\in G$
\begin{equation}
\label{kappa-G}
\pi(h)^* \ka^{Q',hgh^{-1}} = \ka^{Q,g}\,.
\end{equation}
\end{pred}
Let us first mention that if we have the elements $\ka^{Q,g}$
satisfying properties (\ref{kappa}) and (\ref{kappa-G})
the proof of proposition \ref{ONO} is completed.

Indeed, using $\ka^{Q,g}$ we can construct a collection of
morphisms
$$
\psi^{Q,g}\,:\, (\Om^{\bul}(X^Q_g)((\h)),d)\mapsto
(\Om^{\bul}(\bcC^{Q,g}_{2m^Q_g}), D^{Q,g}+\mb)\,,
$$
\begin{equation}
\label{psi-Qg}
\psi^{Q,g} (\eta) =
\eta \, \ka^{Q,g}
\,, \qquad \eta \in \Om^{\bul}(X^Q_g)((\h))
\end{equation}
of complexes such that for any
$h \in G$ the diagram
\begin{equation}
\label{psi-equiv}
\begin{array}{ccc}
  \Om^{\bul}(X^{Q'}_{hgh^{-1}})((\h))  &
\stackrel{\psi^{Q',hgh^{-1}}}{\longrightarrow}&
\Om^{\bul}(\bcC^{Q',hgh^{-1}}_{2m^{Q'}_{hgh^{-1}}}) \\[0.3cm]
\downarrow^{\pi(h)^{\ast}}  & ~  & \downarrow^{\pi(h)^{\ast}}\\[0.3cm]
  \Om^{\bul}(X^Q_g)((\h))  &
\stackrel{\psi^{Q,g}}{\longrightarrow}&
\Om^{\bul}(\bcC^{Q,g}_{2m^Q_g})
\end{array}
\end{equation}
commutes.

Due to (\ref{psi-equiv}) the maps (\ref{psi-Qg})
form a morphism of complexes of $G$-modules
\begin{equation}
\label{psi-tot}
\psi^{tot}\,:\,
\bigoplus_{g,Q} \Om^{2m^Q_g - \bul}(X^Q_g)((\h))
\mapsto
\bigoplus_{q-p= \bul} \bigoplus_{Q,g}
(\Om^{p}(\bcC_{q}^{Q,g}), D^{Q,g}+\mb)\,.
\end{equation}

The fibers of the complex of bundles $\bcC^{Q,g}_{\bul}$
is isomorphic to the normalized Hochschild complex
$\bC(W,Wg)$ of the Weyl algebra with values in the
twisted module $Wg$. Thus, due to proposition \ref{ochen-nado}
and the fact that the evaluations of the elements $\ka^{Q,g}$
are non-trivial cycles in fibers of $\bcC^{Q,g}$,
the map (\ref{psi-tot}) induces an isomorphism
between the $G$-modules
\begin{equation}
\label{psi-iso}
\bigoplus_{Q,g} \Om^{\bul}(X^Q_g)((\h))
\erarrow \bigoplus_{Q,g}
H_{\mb}(\Om^{\bul}(\bcC^{Q,g}))\,,
\end{equation}
where $H_{\mb}$ denotes cohomology of the
double complex $\Om^{p}(\bcC_{q}^{Q,g})$
with respect to the differential $\mb$\,.

Furthermore, proposition \ref{acyc-D} implies that the
map (\ref{la-J}) composed with the projection (\ref{normal})
induces an isomorphism of the $G$-modules
\begin{equation}
\label{la-iso}
J_{\bul} \erarrow \bigoplus_{Q,g}H_{D^{Q,g}}(\Om(\bcC_{\bul}^{Q,g}))\,,
\end{equation}
where $H_{D^{Q,g}}$ denotes cohomology of the
double complex $\Om^{p}(\bcC_{q}^{Q,g})$ with
respect to the differential $D^{Q,g}$\,.

Proposition \ref{acyc-D} also implies that
$\Om^p(X^Q_g, \bcC^{Q,g}_q)$ is acyclic along all the rows
outside the intersection with the zeroth column and
proposition \ref{ochen-nado} implies that
$\Om^p(X^Q_g, \bcC^{Q,g}_q)$ is acyclic along all the
columns outside the intersection with the $2m^Q_g$-th row,
where $2m^Q_g$ is the dimension of $X^Q_g$\,.
Thus, since the subcomplexes $\Om^p(X^Q_g, \bcC^{Q,g}_q)$
of the double complex $\cD_{p,q}$ (\ref{dir-s}) are
placed in the strips $0\le p \le 2m^Q_g$,
we can apply to $\cD_{p,q}$ the standard argument of the
spectral sequence and conclude that the
composition of the morphism of (\ref{q-iso}),
(\ref{la-J}), and (\ref{normal}) is the
first desired quasi-isomorphism in (\ref{dir-sum}) and
the map (\ref{psi-tot}) is the second
desired quasi-isomorphism in (\ref{dir-sum})\,.

~\\
{\bf Proof} {\it of proposition \ref{el-kappa}}~
We set
\begin{equation}
\label{ka-0}
\ka^{Q,g}_0 = \ve_{a_1\dots a_{2m^Q_g}}
\Pi( 1 \otimes y^{a_1} \otimes \dots \otimes
y^{a_{2m^Q_g}}) \in
\Om^0(\bcC^{Q,g}_{2m^Q_g})\,,
\end{equation}
where $\ve_{a_1\dots a_{2m^Q_g}}$ are the
components of the Liouville volume form associated with
the symplectic form on $X^Q_g$,  $y^{a}$ are
fiber coordinates in $TX^Q_g$ and $\Pi$ is the projection
(\ref{normal})\,.

By proposition \ref{ochen-nado}
\begin{equation}
\label{ka-mb}
\mb \ka^{Q,g}_0 =0
\end{equation}
and the evaluations of $\ka^{Q,g}_0$ represent
non-trivial cycles in fibers of $\bcC^{Q,g}$\,.

Unfortunately, $D^{Q,g}\ka^{Q,g}_0$ could be non-zero.
However, due to (\ref{ka-mb}) the element
$\nu_1= D^{Q,g}\ka^{Q,g}_0$ is a $\mb$-cycle
$\mb \nu_1 =0$. Let us show that $\nu_1= D^{Q,g}\ka^{Q,g}_0$
is a $\mb$-boundary.

By definition (\ref{F-conn-Qg})
\begin{equation}
\label{nuu-11}
\nu_1 = \n^{Q,g} \ka_0^{Q,g}  + \frac1{\h}[\cA^{Q,g}, \ka_0^{Q,g}]\,,
\end{equation}
where the commutator is defined
componentwise.

The first term in (\ref{nuu-11}) vanishes since
the symplectic form and hence the Liouville tensor
on $X^Q_g$ are preserved by the connection $\n^{Q,g}$\,.
Therefore,
\begin{equation}
\label{nuu-1}
\nu_1 = \frac1{\h}\ve_{a_1\dots a_{2m^Q_g}}
\sum_{k=1}^{2m^Q_g}
\Pi( 1 \otimes y^{a_1} \otimes\dots \otimes
y^{a_{k-1}}\otimes  [\cA^{Q,g}, y^{a_k}] \otimes
y^{a_{k+1}}\dots \otimes y^{a_{2m^Q_g}})\,.
\end{equation}
But the right hand side of (\ref{nuu-1}) is
the $\mb$-boundary of the following
element\footnote{This fact has a simple
explanation in terms of the action of Hochschild
cochains on Hochschild chains of any associative
algebra \cite{D-Ts}.}:
$$
\nu_2 =- \frac1{\h}\ve_{a_1\dots a_{2m^Q_g}}
\Pi(1 \otimes \cA^{Q,g} \otimes y^{a_1} \otimes \dots \otimes
y^{a_{2m^Q_g}}) -
$$
$$
\frac1{\h}\ve_{a_1\dots a_{2m^Q_g}} \sum_{k=1}^{2m^Q_g}
(-)^{k} \Pi
(1 \otimes y^{a_1}\otimes\dots y^{a_k}\otimes \cA^{Q,g} \otimes
y^{a_{k+1}} \otimes \dots \otimes
y^{a_{2m^Q_g}})\in \Om^1(\bcC^{Q,g}_{2m^Q_g+1})\,.
$$
Thus,
$$
D^{Q,g} \ka^{Q,g}_0 = \mb \nu_2\,.
$$
Due to this observation and the fact that the double
complex $\Om^{p}(\bcC^{Q,g}_q)$ is acyclic along all the
columns outside the intersection with the $2m^Q_g$-th
row we can extend the element (\ref{ka-0}) to the
sum
\begin{equation}
\label{ka-sum}
\tka^{Q,g} = \ka^{Q,g}_0+ \sum_{i=1}^{2m^Q_g}
\ka^{Q,g}_i\,, \qquad
\ka^{Q,g}_i \in \Om^{i}(\bcC^{Q,g}_{2m^Q_g+i})\,,
\end{equation}
such that
$$
(D^{Q,g}+ \mb) \tka^{Q,g}=0\,.
$$

Due to acyclicity of the double
complex $\Om^{p}(\bcC^{Q,g}_q)$ along all the
rows outside the intersection with the zeroth column
we can construct a cycle $\ka^{Q,g}$ which is
placed in the subgroup $\Om^0(\bcC^{Q,g}_{2m^Q_g})$
and is homologically equivalent to (\ref{ka-sum}).

It is clear from the construction that
the difference between $\ka^{Q,g}-\ka^{Q,g}_0$ is
a $\mb$-boundary. Hence the evaluation $(\ka^{Q,g})\Big|_{p}$ at
any point $p\in X^Q_g$ is a non-trivial
$\mb$-cycle in the fiber $(\bcC^{Q,g})_{p}$\,.

It is also obvious that for any $h\in G$
\begin{equation}
\label{kappa0-G}
\pi(h)^* \ka^{Q',hgh^{-1}}_0 = \ka^{Q,g}_0\,.
\end{equation}
Therefore, since we use homotopy operators
which are compatible with the action the general
linear group in fiber coordinates $y^a$,
the constructed elements $\ka^{Q,g}$ satisfy
equation (\ref{kappa-G}) for any $h\in G$\,.

We finish the proof of proposition \ref{el-kappa} and
the proof of theorem \ref{ONA}. $\Box$

\section{Applications to deformation theory}

Let $A$ and $G$ be as in Theorem 1.

\begin{conjecture}\label{conj}
The deformations of the algebras $A[G]$ and $A^G$ are
unobstructed. Thus there exists universal
formal deformations $H_c$, $H_c'$ of these algebras parametrized by
$c\in H^2(A[G])=H^2(A^G)$.
\end{conjecture}

Note that according to Theorem 1,
$$
H^2(A^G)=H^2_{CR}(X/G)=
(H^2(X)\bigoplus \oplus_{g\in G}\oplus_{Q: {\rm codim}X_g^Q=2}
H^0(X_g^Q))^G.
$$
(where the coefficients on the RHS are $\Bbb C((\hbar))$).
Thus the conjecture implies that if $G$ contains symplectic
reflections (i.e., elements $g$ for which $X_g$ contains
components of codimension 2 in $X$), then there exist
``interesting'' deformations of $A^G$, i.e., ones not coming from
$G$-invariant deformations of $A$.

Let us give a few examples motivating this conjecture.

\begin{example} If $G$ is trivial, the conjecture is true,
as follows e.g. from theorem $5.3.3$ in \cite{F1}.
\end{example}

\begin{example}
Let $X=T^*Y$, where $Y$ is a smooth affine variety,
and suppose $G$ acts on $Y$. In this case the conjecture is true
as follows from Theorem 2.13 of \cite{Pasha}, and
$H_c$ is called the Cherednik algebra of $X,G$.
If $X$ is a vector space and $G$ acts linearly, then
$H_c$ is the rational Cherednik algebra from \cite{EG}.
\end{example}

\begin{example} If $X=V$ is a symplectic vector space and $G$
acts linearly then the conjecture holds by Theorem 2.16 of
\cite{EG}, and $H_c$ is the symplectic reflection algebra,
the main object of study in \cite{EG}.
\end{example}

\begin{example} Let $X=V/L$, where $V$ is a symplectic vector
space and $L$ a lattice in $V$ (i.e., $L$ is the abelian group
generated by a basis of $V$). Thus $X$ is an algebraic
torus with a symplectic form. We assume that the symplectic form
is integral and unimodular on $L$. Let $G\subset Sp(L)$ be a
finite subgroup; then $G$ acts naturally on $X$.
In this case we expect that $H_c$ is an
``orbifold Hecke algebra'' defined in \cite{Pasha},
Example 3.11.

More precisely, let $V_{\Bbb R}$
be the tensor product $L\otimes_{\Bbb Z} {\Bbb R}$.
This space has a real symplectic structure $\omega$
obtained by restricting the symplectic
structure on $V$. Pick a $G$-invariant K\"ahler structure on
$V_\Bbb R$ whose imaginary part is $\omega$.
This makes $V_{\Bbb R}$ into a complex vector space,
and thus $X_{\Bbb R}=V_{\Bbb R}/L$ is a complex torus (analytic
abelian variety).
The construction of paper \cite{Pasha} attaches to $X_{\Bbb R},G$
the Hecke algebra ${\mathcal H}_\tau(X_{\Bbb R},G)$,
parametrized exactly by the space $H^2(A[G])$.
We expect this algebra to coincide with $H_c$.

Note that this class of examples includes Cherednik's double
affine Hecke algebras which are attached to affine Weyl groups \cite{Che}.
\end{example}

Motivated by these examples, we propose to call
$H_c$ the symplectic reflection algebra of $X,G$.

\begin{teo}
Assume that $H^3(X,\Bbb C)^G=0$, and for any
$g,Q$ such that $X_g^Q$ has codimension 2, one has
$H^1(X_g^Q,\Bbb C)^{Z(g)}=0$, where $Z(g)$ is the centralizer of
$g$ in $G$. Then conjecture \ref{conj} is true.
\end{teo}

{\bf Proof.}
According to Theorem 1, the assumptions imply that
$H^3(A[G])=H^3(A^G)=0$, so the statement follows from classical
deformation theory. $\Box$

\vskip .05in

Note that in the case when $X=Y^n$, where $Y$ is a smooth
affine symplectic variety, and $G=S_n$, this theorem
was proved in \cite{EO}. In this case, interesting deformations
of $A^G$ arise if and only if $Y$ is an algebraic surface.


\begin{thebibliography}{99}

\bibitem{Alev} J. Alev, M.A. Farinati, T. Lambre,
et A.L. Solotar, Holomologie des invariant d'une
alg$\acute{{\rm e}}$bre de Weyl sous l'action
d'un groupe fini, J. of Algebra, {\bf 232} (2000)
564-577.

\bibitem{Bayen}  F. Bayen, M. Flato, C. Fronsdal, A. Lichnerowicz, and D.
Sternheimer,  Deformation theory and quantization. I.
Deformations of symplectic structures, Ann. Phys. (N.Y.) {\bf 111} (1978), 61;\\
Deformation theory and quantization. II. Physical applications,
Ann. Phys. (N.Y.) {\bf 110} (1978), 111.


\bibitem{Ber}  F.A. Berezin, Quantization, Izv. Akad. Nauk. {\bf 38}
(1974) 1116-1175;\\
General concept of quantization,
Commun. Math. Phys. {\bf 40} (1975) 153-174.


\bibitem{J-L} J.-L. Brylinski, Algebras associated with group actions
and their homology, Brown University preprint, 1987.


\bibitem{J-L1} J.-L. Brylinski, Cyclic homology and equivariant theories,
Ann. Inst. Fourier {\bf 37}, 4  (1987) 15--28.


\bibitem{Br} J.-L. Brylinski, A differential complex for Poisson manifolds,
J. Diff. Geom.  {\bf 28}, 1  (1988) 93--114.


\bibitem{CR} W. Chen and Y. Ruan,
A New Cohomology Theory for Orbifold,
Commun. Math. Phys. {\bf 248} (2004) 1--31,
math.AG/0004129.

\bibitem{Che} I. Cherednik,  Double affine Hecke algebras,
Knizhnik-Zamolodchikov equations, and Macdonald operators,
IMRN {\bf 9} (1992) 171--180.


\bibitem{D-Ts} Yu. Daletski and B. Tsygan, Operations
on cyclic and Hochschild complexes, Methods Func. Anal.
Topology, {\bf 5}, 4 (1999) 62-86.

\bibitem{chains} V. Dolgushev, A formality theorem
for Hochschild chains, to appear in Adv. Math.;
math.QA/0402248.

\bibitem{Hoch-DR} V. Dolgushev, Hochschild Cohomology versus De
Rham Cohomology without Formality Theorems, math.QA/0405177.

\bibitem{DMVV} R. Dijkgraaf, G. Moore, E. Verlinde, and
H. Verlinde, Elliptic Genera of Symmetric Products and Second
Quantized Strings, Commun. Math. Phys.
{\bf 185} (1997) 197--209; hep-th/9608096.


\bibitem{Pasha} P. Etingof, Cherednik and Hecke algebras
of varieties with a finite group action, math.QA/0406499.

\bibitem{EG} P. Etingof, V.Ginzburg, Symplectic reflection algebras,
Calogero-Moser systems, and a deformed Harish-Chandra
isomorphism, Inventiones Math, vol. 147, (2002), p. 243-348

\bibitem{EO} P. Etingof and A. Oblomkov,
Quantization, orbifold cohomology, and Cherednik algebras,
math.QA/0311005.

\bibitem{F} B.V. Fedosov, A simple geometrical construction of deformation
quantization, J. Diff. Geom. {\bf 40} (1994) 213--238.

\bibitem{F1}  B.V. Fedosov, Defomation quantization and index
theory. Akademia Verlag, Berlin, 1996.

\bibitem{FT} B.L. Feigin and B.L. Tsygan,
Cyclic homology of algebras with quadratic relations,
universal enveloping algebras and group algebras,
Lect. Notes in Math. {\bf 1289}, Springer, Berlin, 1987

\bibitem{chiral} E. Frenkel and M. Szczesny,
Chiral de Rham Complex and Orbifolds,
math.AG/0307181.

\bibitem{Nora} N. Ganter, Orbifold genera, product formulas
and power operations, PhD thesis (MIT); math.AT/0407021.

\bibitem{GJ} E. Getzler and J. D. S.  Jones,
The cyclic homology of crossed product algebras,
J. Reine Angew. Math.  {\bf 445}  (1993) 161--174.

\bibitem{GJ1} J. Block, E. Getzler, and J. D. S.  Jones,
The cyclic homology of crossed product algebras. II.
Topological algebras.  J. Reine Angew. Math.  {\bf 466}
(1995) 19--25.

\bibitem{GK}  V. Ginzburg and D. Kaledin,
Poisson deformations of symplectic quotient singularities,
Adv. Math.  {\bf 186}, 1  (2004) 1--57;
math.AG/0212279.


\bibitem{HKR} G. Hochschild, B. Kostant, and A. Rosenberg,
Differential forms on regular affine algebras, Trans. Amer. Math.
Soc. {\bf 102} (1962) 383-408.

\bibitem{Lorenz} M. Lorenz, On the homology of graded algebras,
Commun. Alg. {\bf 20}, 2 (1992) 489-507.


\bibitem{Loday} J.- L. Loday, Cyclic Homology,
Grundlehren der mathematischen Wissenschaften, 301.
Springer-Verlag, Berlin, 1992.

\bibitem{Mont} O.M. Neroslavskii and A.E. Zaleskii,
On simple Noetherian rings, Proc. Akad. Bylorussia S.S.R.
(1975) 38-42. (Russian)


\bibitem{NT} R. Nest and B Tsygan, Algebraic index theorem,
Commun. Math. Phys. {\bf 172}, 2 (1995) 223-262.

\bibitem{Tang} X. Tang, Quantization of non-commutative
Poisson manifolds, PhD thesis (UC Berkeley).


\bibitem{VB} M. Van Den Bergh, A Relation between Hochschild
Homology and Cohomology for Gorenstein Rings, Proc. Amer. Math.
Soc. {\bf 126}, 5 (1998) 1345-1348;\\
Erratum to  ``A Relation between Hochschild
Homology and Cohomology for Gorenstein Rings'',
Proc. Amer. Math. Soc. {\bf 130}, 9 (2002) 2809-2810.


\bibitem{Xu} P. Xu, Fedosov *-products and quantum momentum maps,
Commun. Math. Phys.  {\bf 197}, 1  (1998) 167--197;
q-alg/9608006.

\end{thebibliography}
\end{document}